\documentclass[english]{amsart}

\usepackage{babel}
\usepackage{amstext}
\usepackage{amsmath}
\usepackage{amsfonts}
\usepackage{latexsym}
\usepackage{ifthen}
\usepackage{xypic}
\xyoption{all}
\newcommand{\Pic}{{\rm Pic}}

\newtheorem{lemma1}[equation]{}

\newenvironment{lemma}{\begin{lemma1}{\bf Lemma.}}{\end{lemma1}}

\newenvironment{theorem}{\begin{lemma1}{\bf Theorem.}}{\end{lemma1}}

\newenvironment{proposition}{\begin{lemma1}{\bf Proposition.}}{\end{lemma1}}
\newenvironment{corollary}{\begin{lemma1}{\bf Corollary.}}{\end{lemma1}}

\newenvironment{remark}{\begin{lemma1}{\bf Remark.}\rm}{\end{lemma1}}
\newenvironment{definition}{\begin{lemma1}{\bf Definition.}}{\end{lemma1}}

\newenvironment{conjecture}{\begin{lemma1}{\bf Conjecture.}}{\end{lemma1}}

\makeatletter
\ifnum\@ptsize=0  \fi \ifnum\@ptsize=2  \fi 

\newcommand\sE{{\mathcal E}}
\newcommand\sF{{\mathcal F}}

\newcommand\sO{{\mathcal O}}

\newcommand\bR{{\mathbb R}}

\newcommand\bQ{{\mathbb Q}}
\newcommand\bN{{\mathbb N}}

\newcommand\bP{{\mathbb P}}

\begin{document}
\title {Strictly nef divisors} \author{Fr\'ed\'eric Campana, Jungkai A. Chen and Thomas Peternell}
\date{october 25, 2005}

\maketitle

\vskip .2cm \noindent

\vspace*{-0.5in}\section*{Introduction}
Given a line bundle $L$ on a projective manifold $X$, the Nakai-Moishezon criterion says that $L$ is ample if and only if
$$ L^s \cdot Y > 0 $$
for all $s$ and all irreducible subvarieties $Y \subset X$ of
dimension $s.$ Examples show that it is not sufficient to assume
that $L \cdot C > 0$ for all curves; line bundles with this
property are called {\it strictly nef }. If however $L = K_X$ is
strictly nef, then standard conjectures predict that $K_X$ is
already ample; this is proved by ``Abundance'' in dimension up to
3 (Kawamata, Miyaoka). If $L = -K_X$ is strictly nef in dimension 3, then Serrano
[Se95] showed that $-K_X$ is ample, i.e. $X$ is a Fano threefold. This lead
him to set up the following
\begin{conjecture} Let $X_n$ be a projective manifold and $L$ a strictly nef line bundle on $X$. Then $K_X + tL $ is ample for any real
$t > n+1.$
\end{conjecture}
Serrano established the conjecture in dimension 2, and also in dimension 3, with the following possible exceptions:
\begin{itemize}
\item $X$ is Calabi-Yau and $L \cdot c_2(X) = 0;$
\item $X$ is uniruled with irregularity $q(X) \leq 1$, in particular $X$ is rationally connected;
\item $X$ is uniruled with irregularity $q(X) = 2$ and $\chi(\sO_X) = 0$.
\end{itemize}

As said, he also settled the case $L = -K_X$ in dimension 3.

In this paper we rule out the two last cases and establish also results in higher dimensions:

\begin{theorem} Let $X_n$ be a projective manifold and $L$ a strictly nef line bundle on $X$.
Then $K_X + tL $ is ample if $t > n+1$ in the following cases.
\begin{enumerate}
\item $\dim X = 3$ unless (possibly) $X$ is Calabi-Yau with $L \cdot c_2 = 0;$
\item $\kappa (X) \geq n-2;$
\item $\dim \alpha(X) \geq n-2$, with $\alpha: X \to A$ the Albanese map.
\end{enumerate}
\end{theorem}

Statement 2) (resp. 3), resp. 1)) will be proved in \S 2 (resp. \S 3, resp. \S4-5).

The remaining three-dimensional case that $X$ is Calabi-Yau with $L \cdot c_2 = 0$ is a very hard problem in Calabi-Yau theory and definitely requires
very different methods.

\tableofcontents

\section{Basic definitions, known results and main problems}

For technical reasons we have to consider not only strictly nef line bundles, but
also a slight generalization of this notion.
\setcounter{equation}{0}
\begin{definition} Let $X$ be a normal projective variety.
\begin{enumerate}
\item
A line bundle $L$ over $X$ is strictly nef, if $L \cdot C > 0$ for
all irreducible curves $C \subset X.$
\item $L$ is almost strictly nef, if
there is a normal projective variety $X'$, a surjective birational
holomorphic map $f: X \to X'$ and a strictly nef line bundle $L'$ on
$X'$ such that $L = f^*(L').$
\end{enumerate}
\end{definition}

The main problem on strictly nef line bundles is Serrano's:

\begin{conjecture} Let $X_n$ be a projective manifold and let $L$ be a strictly nef line bundle on $X$.
Then $K_X + tL$ is ample for $t > n+1$.
\end{conjecture}

\

{\bf Remark:}  {\it More generally, one might conjecture that if $X_n$ is a normal projective
variety with canonical singularities and index $i(X)$, and if  $L$ is a strictly nef
line bundle on $X$, then $K_X+ tL$ is ample for all $t > i(X)(n+1).$

By definition, the index $i(X)$ is the smallest number $i$ such
that $iK_X$ is Cartier. One could add (in the smooth case)  that
$K_X + nL$ is always nef, and not ample if and only if $X = \bP_n,
L = \sO_X(1).$}

\

It is known since a long time that strictly nef divisors need not be ample; even if
moreover big. See Ramanujam's example in [Ha70].

There are however three important special cases of the conjecture,
namely when $L = K_X$ (resp. $L = -K_X)$, resp. $K_X \equiv 0$. In
the first case the abundance conjecture predicts that $mK_X$ is
spanned for a suitable large $m$ so that $K_X$ will be ample as
soon as $K_X$ is strictly nef. This is known in dimension up to
$3.$ In the second case $X$ should be Fano if $-K_X$ is strictly
nef. In the last case, $L$ should be ample.

\begin{remark}  Perhaps the best justification for the above
conjecture (1.2) is that it holds for $L$ if and only if

$$ L^{\perp} \cap K_X^{\perp} \cap {\overline {NE}}(X) = \{0\},$$
in $N_1(X)$, see Proposition 1.4 below. So the conjecture should
be viewed as a statement on the cone ${\overline {NE}}(X),$ at the
points where the intersection number with $K$ and $L$
simultaneously vanish. Observe thus that the crucial cases are
precisely the three ``special" cases above, where $L= K_X$,
$L=-K_X$, and $ K_X \equiv 0$.

Notice also that, if ${\overline {NE}}(X)$ is generated by
the classes of irreducible curves (i.e. without taking limits), then
the conjecture is
true since $K_X+t L$ is again strictly nef, for $t > (n+1)$ (1.6).
This holds in particular if $X$ is Fano.
\end{remark}

By $ME(X)$ we will always denote the cone of movable curves. Its closure is the cone dual to the
cone of effective divisors;
see [BDPP04] for details.

\begin{proposition} Let $L$ be strictly nef and $\alpha \in {\overline {NE}}(X) $
such that $(K_X + tL) \cdot \alpha = 0$ ($ t > n+1)$. Then
\begin{enumerate}
\item $K_X \cdot \alpha = L \cdot \alpha = 0.$
\item $\alpha \in \partial ME(X)$ for a suitable choice of $\alpha.$
\end{enumerate}
\end{proposition}

\begin{proof} (1) Suppose $L \cdot \alpha \ne 0.$ Then $L \cdot \alpha > 0$ and $K_X \cdot \alpha < 0.$
By the cone theorem we can write
$$ \alpha = \sum_{i=1}^N a_i C_i + R$$
with $C_i$ extremal and $K_X \cdot R \geq 0.$
Since $-K_X \cdot C_i \leq n+1$, and $tL.C_i\geq t>n+1$, for all $i$,  we have $(K_X + tL) \cdot C_i > 0,$ which gives a contradiction. \\
(2) If there is no nonzero $  \alpha \in \partial ME(X)$ with $(K_X + tL) \cdot \alpha = 0,$ then by
[BDPP04] $K_X + tL$ is big. But then $K_X + tL$ is ample, by (1.6(2)) below.
\end{proof}

The following cases have been settled by Serrano [Se95]

\begin{theorem} ({\it Serrano}) \begin{enumerate}
\item Let $X$ be a irreducible reduced projective Gorenstein surface
and $L$ strictly nef on $X$. Then
$K_X+tL$ is ample for any real $t > 3.$
\item Let $X$ be a smooth projective threefold and $L$ strictly nef.
Then $K_X+tL$ is ample for
$t > 4$ with the
following possible exceptions only: $X$ is Calabi-Yau and $L \cdot c_2 =
0;$ or $X$ is uniruled with $q \leq 1$; or
$X$ is uniruled, $q =  2$ and $\chi(\sO_X) = 0.$
\end{enumerate}
Moreover $X$ is Fano as soon as $-K_X$ is strictly nef.
\end{theorem}

The following more technical results are also due to Serrano.

\begin{proposition} Let $X$ be an n-dimensional connected projective
manifold and $L$ a strictly nef line bundle on $X$.
\begin{enumerate}
\item  For every real number $t > n+1,$ $K_X+tL$ is a strictly nef
$\bR-$divisor. This also holds for $t >> 0$ and $X$ a normal projective variety with only canonical singularities.
\item If $K_X+tL$ is not ample for some real number $t > n+1,$ then
$K_X^j \cdot L^{n-j} = 0$
for all $j \geq 0.$ So
if $(K_X+tL)^n \ne 0$ for some real number $t > n+1$ (i.e. if $ K_X+tL$
is big and strictly nef),
then $K_X+tL$ is ample.
\item If $\dim X = 3$ and $\vert pK_X+qL\vert$ contains an effective non-zero divisor for some integers $p,q$,
then $K_X+t L$ is ample for $t > 4.$
\end{enumerate}
\end{proposition}

The last proposition says in particular that to prove
Conjecture 1.2 for each $t > n+1,$ it is sufficient to prove it for
some
positive integer $t> n+1.$

\section{Results in case of positive Kodaira dimension}
\setcounter{equation}{0}
If $X$ is of general type, then Conjecture (1.2) easily holds:

\begin{proposition} Let $X$ be a projective n-dimensional manifold
with $\kappa (X) = n.$ Let $L$ be strictly nef on $X$. Then
$K_X+(n+1)L$ is ample.
\end{proposition}

\begin{proof} Let $t > n+1$ be a rational number. By (1.6), $K_X+tL$ is
strictly nef. Then
$2(K_X+tL) - K_X$ is big and nef, hence by the base point free
theorem, $K_X+tL$ is semi-ample and strictly nef, hence ample.
\end{proof}

If $X$ is not of general type, things are more complicated; here we
want to use the Iitaka fibration.
For technical reasons we introduce

\begin{conjecture} ${\bf (C_d):}$ Let $F_d$ be a projective manifold
with $\kappa (F) = 0. $ Let $L$ be almost strictly nef.
Then, $K_F+tL$ is big for $t > d+1.$
\end{conjecture}

\begin{theorem} Let $X$ be an n-dimensional connected projective
manifold with $\kappa (X) = k \geq 0.$
Let $L$ be a strictly nef divisor on $X$. Suppose that
$C_d$ holds for $d = n-k.$ Then, $K_X+tL$ is ample for $t > n+1.$
\end{theorem}

\begin{proof} Let $f: X \rightharpoonup Y$ be the Iitaka fibration; we may
assume $\dim Y = k \geq 1$, because
otherwise there is nothing to prove.
Let $\pi: \hat X \to X$ be a sequence of blow-ups such that the
induced map $\hat f: \hat X \to Y$ is holomorphic, and moreover
we can write:
$$ \pi^*(mK_X) = \hat f^*(A) + E  \eqno (*)$$
with an ample divisor $A$ on $Y$ and an effective divisor $E.$
We also have an equality: $K_{\hat X} = \pi^*(K_X) + E'$, for some
effective $E'$. Let us set: $\hat L = \pi^*(L).$ \\
By $(C_d)$ applied to the general fiber $F$ of $\hat f$, the divisor
$K_F + t \hat L$ is big, for $t > d+1,$ even. \\
Thus $\pi^*(K_X) + E' + t \hat L$ is $\hat f-$big. The following
Lemma (2.4) therefore applies,
with $N = \pi^*(K_X) + t \hat L$
and with $D = E'$ and shows the bigness of
$$ B:=\pi^*(K_X) + E' + t \hat L + \hat f^*(A)=(f^*(A)+E)+(\pi^*(K_X)+E'+t\hat{L})$$
Thus by (*), $\pi^*((m+1)(K_X)) + E' + t \hat L = B+E $ is big,
and so does
$$ \pi^*((m+1)K_X) + E' + (m+1)t \hat L=(B+E)+mt \hat L,$$
being the sum
of two divisors, $B+E$, which is big, and
$mt \hat L$, which is nef. \\
Therefore $K_{ \hat X} + t \hat L$ is also big and thus ample, by
(1.6).
\end{proof}

\begin{lemma} Let $g: X \to Y$ be a holomorphic map of projective
varieties. Let $A,N,D$ be $\bQ$-divisors,
with $A$ ample
on $Y$, $N$ nef on $X$ and $D$ effective on $X$. Suppose that $D+N$
is $g-$big, i.e. big on the
general fiber. Then $D+N+g^*(A)$ is big.
\end{lemma}

\begin{proof} Choose and fix $k$ large such that $D+N+g^*(kA)$ is big. \\

(This is a standard fact, seen as follows, by a relative version of
Kodaira's Lemma: let $H$ be $g-$ample on $X$. Then choose $m$ such
that $g_*(m(D+N)-H)$ has positive rank.
This is obviously possible, by the coherence of direct image sheaves,
since $(D+N)$ is $g-$big. See [KMM87,0-3-4], for example.
Now choose $k$ large enough, such that $g_*(m(D+N)-H) + kA$
has a section. Thus $E:=m(D+N)-H + g^*(kA)$ is effective, and $m(D+N)
+ g^*(kA) =H+E$ is of
the form: ample plus effective, and thus big, as claimed).
\vskip .2cm \noindent
Thus $D+N+g^*(kA) = aH+E$ with $H$ ample, $a$ a positive rational
number and $E$ an effective
$\bQ-$divisor. Since $N$ is nef, $N+\epsilon H$ is ample for all
positive numbers $\epsilon;$ choose
$\epsilon$ such that $(k-1)\epsilon < a.$

Next observe, introducing the effective divisor $E' =
(k-1)(D+N+\epsilon H)$ that
$$(k-1)(D+N+\epsilon H) + D + N + g^*(kA) = E' + aH + E =: aH + E''$$
with $E''$ effective.
On the other hand,
$$(k-1) (D+N+\epsilon H) + D + N + g^*(kA) = k(D+N+g^*(A)) +
(k-1)\epsilon H, $$
hence, substracting $(k-1).\epsilon.H$ from both sides, we get the equality:
$$ k(D+N+g^*(A))= aH + E'' - (k-1)\epsilon H = (a-(k-1)\epsilon)H + E''.$$
Since $(a-(k-1)\epsilon) > 0$, by the choice of $\epsilon$, the right
hand hand side divisor is big,
hence $D+N+g^*(A)$ is also big.
\end{proof}

Conjecture $C_1$ being obvious, we are now going to prove $C_2.$

\begin{proposition} Let $X$ be a smooth projective surface with
$\kappa (F) = 0$, $L$ almost strictly nef.
Then $K_X+tL$ is big for $t > 3.$
\end{proposition}

\begin{proof} Fix a rational number $t > 3$ and suppose that $K_X + tL$ is not big.

By blowing down the $(-1)-$curves $E_i$ with $L \cdot E_i = 0,$ we
may assume that
$K_X + tL$ is nef. So if $(K_X + tL)$ is not big, we must have $(K_X
+ tL)^2 = 0.$
If $L^2 > 0,$ then of course our claim is clear, so suppose $L^2  = 0.$ Hence
$$( K_X^2 + 2t K_X) \cdot L = 0.$$
This holds also for all rational numbers $3 < t_0 < t,$ because
otherwise $K_X + t_0L$ would
be big and then also $K_X + tL$ is big.

Thus: $K_X^2=0=K_X\cdot L$. The surface $X$ is thus minimal.
By taking a finite \'etale cover, we can assume $X$ to be either an
abelian, or a $K3$-surface.

But the argument used in [Se95] for abelian varieties shows that an almost
strictly nef divisor on an abelian variety is ample.

On the other hand, Riemann-Roch shows that a nef line bundle on a
$K3$-surface is either
effective or trivial. An effective almost strictly nef line bundle on
a surface is immediately seen to be big, and thus ample.
\end{proof}

\begin{remark} Claim $(C_2)$ trivially holds also on surfaces of
general type and is
very easily checked in case $\kappa = 1.$ It should also hold in case
$\kappa = - \infty$ but
we don't need it.
\end{remark}

Combining (2.3) and (2.5) we obtain:

\begin{corollary} Let $X$ be an n-dimensional connected projective
manifold with $\kappa (X) \geq n-2.$ Let $L$ be
a strictly nef line bundle on $X$. Then $K_X+t L$ is ample for $t > n+1.$
\end{corollary}

\section{The Albanese map}
\setcounter{equation}{0}

We now study Conjecture 1.2 on projective manifolds $X$ with $q(X)
> 0.$ Since our most complete result is in dimension 3,
we will do this case first and then examine what can be done in
higher dimensions.

\begin{theorem} Let $X$ be a smooth projective threefold, $L$ strictly nef.
Suppose there exists a non-constant map $g: X \to A$ to a abelian variety.
Then
$K_X + tL$ is ample for $t > 4.$
\end{theorem}

\begin{proof} Let $D_t:= 2K_X+2tL$. We claim that $\mathcal{F} = \mathcal{F}_t:=g_*(2K_X+2tL)$
satisfies a {\it generic vanishing theorem} (cf. \cite{Ha04}
Theorem 1.2) for $t$ a sufficiently large integer. That is, we have a chain of inclusions
$$ V^0( \mathcal{F}) \supset V^1(\mathcal{F})\ldots \supset
V^n(\mathcal{F}),$$ where
$$V^i(\mathcal{F}):=\{ P \in {\rm
Pic}^0(X)| h^i(A, \mathcal{F} \otimes P) \ne 0\}.$$

Grant the claim for the time being. Since $\mathcal{F}$ is a
non-zero sheaf for $t \gg 0$, one concludes that $V^0(\mathcal{F})
\ne \emptyset$. For otherwise, $V^i(\mathcal{F}) =\emptyset$ for
all $i$, which implies that the Fourier-Mukai transform of $\mathcal{F}$
is zero. This is absurd e.g. by [Mu81,2.2].

Therefore $h^0(X, 2K_X+2tL+P) \ne 0$ for some $P \in \Pic^0(X)$. Now
$L'$ be a divisor such that $2tL'=2tL+P$, then by Proposition
1.6.(3), $K_X+tL'$ is ample for $t > 4$ and hence so is
$K_X+tL$ (notice that if $2K_X + 2tL + P$ has a section without zeroes, then $-K_X $ is strictly nef,
hence $X$ is Fano and $q(X) = 0,$ so that 1.6(3) really applies).

To see the claim, first note that $K_X+t_0L$ is $g$-big for some
$t_0 >0$ (1.5(1)). Fix any ample line bundle $H$ on $A$. By Lemma 2.4, it is
easy to see that $a(K_X+2tL)+g^*H$ is nef and big for $a>0$ and $t
>t_0$. Set $D_0:= 2(K_X+2tL)+g^*H$. then $D_0-K_X$ is again
nef and big. By the Base Point Free Theorem, $mD_0$ is spanned
for some $m \gg 0$. Take $D$ a general smooth member in $|mD_0|$.
Then we have $$ 2K_X+2tL+g^*H \equiv
K_X+\frac{2m-1}{2m}g^*H+\frac{1}{2m}D,$$ where $(X,\frac{1}{2m}D)$ is
klt. By the vanishing theorem of Koll\'ar, we have
$$ H^j(A, g_*(2K_X+2tL) \otimes H))=0, \text{ for all } j >0$$
and moreover
$$ H^j(A,\mathcal {F} \otimes H \otimes P) = H^j(A, g_*(2K_X+2tL) \otimes H \otimes P))=0 \text{ for all } j >0, \  P \in \rm{Pic}^0(X).$$
In other words, {\it per definitionem} the sheaf $\mathcal{F} \otimes H$ is $IT^0$ for all
ample line bundles $H$.

Next, let $M$ be any ample line bundle on the dual abelian variety
$\hat{A}$ and $\phi: \hat{A} \to A$ is the isogeny defined by $M$.
Let $\hat{M}$ be the Fourier-Mukai transform of $M$ on $A$ and let
$\hat{M}^\vee$ be its dual. By \cite{Mu81} Proposition 3.11,
$$ \phi^*( \hat{M}^\vee) \cong \oplus^{h^0(M)} M.$$

Let $\hat{g}: \hat{X}:=X \times_A \hat{A} \to \hat{A}$ be the base
change with $\varphi: \hat{X} \to X$ being \'etale. Clearly,
$K_{\hat{X}}=\varphi^*K_X$ and $\varphi^*L$ is strictly nef on
$\hat{X}$. Let $\mathcal{G}:=\hat{g}_*(\varphi^*(2K_X+2tL))$. By
applying the above argument to $\varphi^*D_t$, we see that
$\mathcal{G} \otimes M$ is $IT^0$ for all $M$. Thus

$$
\begin{array}{rl}
 \phi^*(\mathcal{F} \otimes \hat{M}^{\vee}) & =
\phi^*( g_*( D_t \otimes \hat{M}^{\vee}))\\
&= \phi^*g_*( ( D_t \otimes g^* \hat{M}^{\vee})) \\
&= \hat{g}_* \varphi^*( ( D_t \otimes g^* \hat{M}^{\vee})) \\
&= \hat{g}_*  (\varphi^* D_t \otimes \varphi^*g^*
\hat{M}^{\vee})\\
&= \hat{g}_*  (\varphi^* D_t \otimes \hat{g}^* \phi^*
\hat{M}^{\vee})\\
&=\hat{g}_*  (\varphi^* D_t \otimes \hat{g}^* ( \oplus M))\\
&=\oplus (\hat{g}_* \varphi^* D_t \otimes   M)
\end{array}
$$ which
is $IT^0$.

Since $\mathcal{O}_A$ is a direct summand of $\phi_* \phi^*
\mathcal{O}_A$, it follows that $\mathcal{F} \otimes
\hat{M}^{\vee}$ is $IT^0$. By \cite{Ha04} Theorem 1.2, our claim
follows.
\end{proof}

\begin{corollary}
Let $X$ be a smooth projective threefold, $L$ strictly nef. Suppose that $\tilde q(X) > 0.$
Then $K_X + tL$ is ample for $t > 4.$
\end{corollary}

\begin{proof} By the previous theorem we only have to treat the case that $q(X) = 0$.
Then we choose a finite \'etale cover $h: \tilde X \to X$ such that $q(\tilde X) > 0.$
Hence $K_{\tilde X} + h^*(L) $ is ample for $t > n+1$ and so does $K_X + tL.$
\end{proof}

\begin{remark} There are two obstacles for extending Theorem 3.1 to all dimensions.
The first is the use of 1.6(3) which has to be extended to higher dimensions. We will do this
below.
The second is the $g-$bigness of $K_X + tL.$ This means that $K_F + tL_F$ is big for the
general fiber $F$ of $g.$ Thus we need to argue by induction on the dimension, but of course
we are far from proving the conjecture for arbitrary manifolds (with vanishing irregularity).
\end{remark}

\begin{lemma} Let $X$ be an irreducible reduced projective Gorenstein variety with desingularization
$\pi: \hat X \to X.$
Let $g: X \to A$ non-constant and $L$ be a strictly nef line bundle on $X$. Suppose that $K_X + t_0 L$ is $g-$big for
some $t_0$ and set $\hat L = \pi^*(L).$
\begin{enumerate}
\item The sheaf $\hat \sF = g_*\pi_*(2K_{\hat X} + 2t \hat L) $ satisfies the generic vanishing theorem
$$ V^0(\hat \sF) \supset V^1( \hat \sF) \ldots \supset
V^n(\hat \sF).$$
\item If $t \gg 0,$ then there exists $P \in {\rm Pic}^0(X)$ such that
$$ H^0(2K_X + 2tL + P) \ne 0.$$
\end{enumerate}
\end{lemma}

\begin{proof} (1) This is just what the second part of the proof of Theorem 3.1 gives.
(2) By (1) and the first arguments in the proof of Theorem 3.1 we obtain that
$$ H^0(2K_{\hat X} + 2 \hat L + \hat P) \ne 0 $$
for some $\hat P \in {\rm Pic}^0(\hat X).$ Since $\hat P$ comes from $A$, it is of the form $\pi^*(P).$
Moreover we have $\pi_*(2K_{\hat X}) \subset  2K_X$ since $X$ is Gorenstein, hence claim (2) follows.
\end{proof}

\begin{theorem} Let $X_n$ be an irreducible reduced projective variety with a non-constant map $g: X \to A.$ Let $L$ be a strictly nef line
bundle on $X$ and assume that $K_X + t_0L$ is $g-$big for some $t_0$ (e.g., Conjecture 1.2 holds in dimension $< n).$
Then $K_X +  tL $ is ample for $t > n+1.$
\end{theorem}

\begin{proof} We prove the claim by induction on $n.$
Since we argue numerically, we may ignore $P$ and choose by Lemma
3.4
$$D = \sum m_i D_i \in \vert 2K_X + 2tL \vert $$
for large $t.$ We may select a component, say $D_1,$ such that $\dim g(D_1) \ne 0$ and consider the non-constant
map $g_1: D_1 \to A.$ By induction $K_{D_1} + t L_{D_1}$ is ample, if $t > n.$ We now adopt the methods of [Se95,3.1].
The equation
$$  0 = D \cdot (K_X + tL)^{n-1} $$
leads, via the nefness of  $K_X + tL,$ to
$$ 0 = D_1 \cdot K_X^j \cdot L^{n-1-j}, \ 0 \leq j \leq n-1. \eqno (*)  $$
Choose $k$ such that $L_{D_1}^{n-k} \not \equiv 0,$ e.g., $k = n-1.$ Then
$$ 0 < (K_{D_1} + tL_{D_1})^k \cdot L_{D_1}^{n-k} = D_1 \cdot  (K_{D_1} + tL_{D_1})^k \cdot L_{D_1}^{n-k}. $$
By (*) we obtain
$$ 0 < D_1^2 \cdot L^{n-2}. $$
On the other hand, (*) yields
$$ 0 = D_1 \cdot (2K_X + tL) \cdot L^{n-2} = D_1 \cdot (\sum m_i D_i) \cdot L^{n-2} = $$
$$ = m_1 D_1^2  \cdot L^{n-2} + \sum_{i \geq 2} m_i D_1 \cdot D_i \cdot L^{n-2} \geq  m_1 D_1^2 L^{n-2},$$
a contradiction.
\end{proof}

\begin{corollary} Let $X_n$ be a projective manifold with Albanese map $\alpha: X \to A.$ Let $L$ be strictly nef on $X.$
Assume that $\dim \alpha (X) \geq n-2,$ or $\dim \alpha (X) =  n-3 $ but the general fiber $F$ is not Calabi-Yau with
$L_F \cdot c_2(F) = 0.$ Then $K_X + tL$ is ample for $t > n+1.$
\end{corollary}

\section{Fano fibrations}

\setcounter{equation}{0}

We shall now (in particular) complete the proof of Theorem 0.2 (1). Observe that due to 1.5(2), 2.3 and 3.1, the only cases left are uniruled
threefolds with $q=0$. These cases are thus settled by 4.1-2 and 5.1-2 below.

\

In this section we settle the cases of del Pezzo fibrations over curves and elementary conic bundles over surfaces.

\begin{proposition} Let $X$ be a smooth projective threefold, $L$ strictly nef on $X.$
Suppose that $X$ carries an extremal contraction $f: X \to B$ to a
curve $B.$ Then $K_X + tL $ is ample for large $t.$
\end{proposition}

\begin{proof} Since $K_{X} + tL$ and $L$ are strictly nef
and since $\rho(X) = 2,$ $K_{X} + tL$ is clearly ample for large
$t$ unless $-K_{X}$ and $L$ are proportional. Hence $X$ is Fano by
Serrano's theorem (1.5) which ends the proof.
\end{proof}

\vskip .2cm

\begin{proposition} Let $f: X_3 \to S$ be a conic bundle with $\rho(X/S) = 1.$
If $L$ is strictly nef on $X$, then $K_X + tL$ is ample, for $t > 4.$
\end{proposition}

\begin{proof} By Corollary 3.2 we may assume that $q(S) = 0, $ even after a finite \'etale cover of the smooth surface $S.$ \\
Since $\rho(X/S) = 1,$ we find a positive number $t_0$ such that
$$ K_X + t_0L = \phi^*(K_S+M) $$
with a $\bQ-$divisor $M$ on $S.$
We cube the
equation $t_0L = -K_X + f^*(K_S+M) $ to obtain
$$ 0 = 3K_X^2 \cdot f^*(K_S+M) - 3 K_X \cdot f^*(K_S+M)^2 = 3 K_X^2 \cdot f^*(K_S+M) + 6 (K_S+M)^2.$$

From $K_X \cdot t_0L^2 = 0$ we get $K_X^2 \cdot f^*(K_S+M) = 0,$
hence in total
$$ (K_S+M)^2 = 0.$$
By applying (1.4), we find
 $\alpha \in  {\overline {ME}}(X)$ such that
$$ K_X \cdot \alpha = L \cdot \alpha = 0,$$
in particular $  D \cdot \alpha = 0.$ Introducing $\gamma = \phi_*(\alpha) \in {\overline {ME}}(S),$ we obtain
$$ (K_S + M) \cdot \gamma = 0. $$
Notice that  ${\overline {ME}}(S) $ is nothing than the nef cone, so $\gamma$ is a nef class.
Next notice that we may choose $\gamma$ rational. In fact, since the rational points are dense in the nef cone on $S$ and since neither
$K_S+M$ nor $-(K_S+M)$ are strictly positive functionals on the nef cone, we find rational points $x$ and $y$ in the nef cone such that
$$ (K_S+M) \cdot x \geq 0; \ (K_S+M) \cdot y \leq 0.$$
We may assume strict inequality in both cases, otherwise we are already done. Then choose $\lambda > 0$ such that
$$ (K_S+M) \cdot (x+\lambda y) = 0.$$
Noticing that $\lambda \in \bQ,$ we may substitute $\gamma$ by $x + \lambda y.$ Now multiply $\gamma$ suitably to obtain a nef line
bundle $G$
such that
$$ (K_S+M) \cdot G = 0.$$
If now $G^2 > 0,$ then Hodge Index gives $K_S+M = 0$, so that
$H^0(m(K_X+t_0L)) \ne 0$ for positive integers $m$ such that $mt_0
\in \bN.$

Thus we may assume that $G^2=0$. Together with $(K_S+M)^2=(K_S+M)
\cdot G=0$, one has $(K_S+M+ \tau G)^2=0$ for all $\tau$.

Let $C \subset S$ be an irreducible curve. Then
$$ t_0^2  L^2 \cdot f^*(C) = (f^*(K_S+M) - K_{X})^2 \cdot f^*(C) = -2f^*(K_S+M)\cdot K_{X} \cdot C + K_{X}^2 \cdot f^*(C) = $$
$$ = 4 (K_S+M) \cdot C - (4K_S + \Delta) \cdot C = (M - \Delta) \cdot C. \eqno (1)$$
The last equation is explained as follows. Outside the singular
locus of $S,$ the map $f$ is a conic bundle; let $\Delta$ denote
the closure of the discriminant locus. Then it is well-known that
$$f_*(K_{X}^2) = -(4K_S + \Delta).$$
Now we restrict ourselves to curves $C$ with $C^2 \geq 0.$ Then
clearly $L^2 \cdot f^*(C) \geq 0,$ hence
$$ (M-\Delta) \cdot C \geq 0, \eqno (2)$$
in particular $$(M-\Delta) \cdot G \geq 0. \eqno (3)$$ Moreover we
have a strict inequality in (2) unless $C_0 = \emptyset $ and $L^2
\cdot f^*(C) = 0.$ The inequality (3) says im particular that $M$
is pseudo-effective. Thus the equation $(K_S+M) \cdot G = 0$
forces $\kappa (S) \leq 1.$

\vskip .2cm \noindent
{\bf (I)} We first assume $\kappa (S) = - \infty.$ Then $S$ is a rational surface.
The case that $S=\mathbb{P}^2$ is easy. So we may assume that
$\pi: S=S_n \to S_{n-1} \to ... \to S_0$ is a succession of blow-ups,
where $S_0$ is a ruled surface with minimal section
$C_0$ that $C_0^2=-e$.

Now we write $$K_S+M = \pi^*( \alpha_1 C_0 + \beta_1 F ) + E_1$$
$$G=\pi^*( \alpha_2 C_0 + \beta_2 F ) + E_2,$$
where $E_1,E_2$ are divisors supported on exceptional curves.

If $\alpha_2=0$, then it is clear that $E_2=0$ and $G=\beta_2
\pi^*F$. Then $(K_S+M) \cdot G =0$ gives $\alpha_1=0$ and
$(K_S+M)^2=0$ gives $E_1=0$. So $K_S+M= \beta_1 \pi^* F$, and we
are done.

If $\alpha_2 \ne 0$, take $\tau =\frac{ -\alpha_1}{\alpha_2}$,
then
$$ K_S+M+\tau G = (\beta_1 +\tau \beta_2) \pi^*F + E_1 + \tau
E_2.$$ $(K_S+M+\tau G)^2=0$ gives $(E_1 + \tau E_2)^2=0$. It
implies $E_1 + \tau E_2 =0$ by the negativity of intersection form
of exceptional divisors.

Let $\delta:=\beta_1 +\tau \beta_2$. If $\delta \ne 0$, then
$K_S+M= -\tau G + \delta \pi^* F$. Again, one has $\alpha_2=G
\cdot \pi^*F=0$ which is absurd.

Therefore, $K_S+M = -\tau G$. Now since, by $(3)$, $M \cdot G \ge
0$, we have $K_S \cdot G \le 0$. By Riemann-Roch and the obvious vanishing 
$H^0(K_S-G) = 0,$ we have 
$$h^0(S,G)  \ge \chi(\mathcal{O}_S)=1.$$ Hence $G$ is effective. $G$ is non-zero for
otherwise $K_X + t_0L \equiv 0,$ hence $-K_X$ is strictly nef and thus $X$ is Fano. Therefore $m(K_X+t_0L)$
is effective for some $m \in \mathbb{Z}$ and we are done in Case
(I).

\vskip .2cm \noindent
{\bf (II)} Now suppose that $\kappa (S) \geq 0.$
 Let $$ \sigma: S \to S_0$$ be the
minimal model. Since $\kappa (S) \geq 0,$ we conclude by (2) that
$$ K_S \cdot G = M \cdot G = 0.$$
Hence $$ 0 = \sigma^*(K_{S_0}) \cdot G + \sum a_iA_i \cdot G$$
with $A_i$ the $\sigma-$exceptional curves and $a_i$ suitable positive rational numbers.
Thus $G = \sigma^*(G_0)$ with a nef line bundle $G_0$ on $S_0$; observe that $K_{S_0} \cdot G_0 = 0$ and that $G_0^2 = 0.$

\vskip .2cm \noindent Suppose that $\kappa (S) = 1.$ Then we
consider the Iitaka fibration $g: S_0 \to B$ to the curve $B$
(necessarily $B = \bP_1).$ We conclude that $G_0$ is a sum of
fibers of $g$. Thus $G$ is a sum of fibers of $g \circ \sigma.$
Now consider the composed map $h: X \to B.$ Then it follows that
$h_*(\alpha)$ consists of finitely many points. This means that we
can find a fiber of $h$ such that $K_X+tL \vert F$ is not ample
for large $t$. Thus $K_F+tL_F$ is not ample. If (the reduction of)
$F$ is irreducible, this contradicts (1.6). If $F_i$ is a
component of $F$ with multiplicity $a_i,$ then $a_iK_{F_i} +
tL_{F_i}$ is a subsheaf of $K_F + tL_F \vert L_{F_i}$, and the
contradiction is the same.

\vskip .2cm \noindent
Finally we have to treat the case $\kappa (S) = 0.$ Here we may assume that $S_0$ is K3.
If $G_0^2 = 0$, then by Riemann-Roch $\kappa (G_0) = 1.$ Hence some multiple of $G_0$ is spanned, defining a morphism $g: S \to B.$
Since the divisor $M_0$ must be
supported on fibers of $g,$ so does $\Delta.$ Thus we conclude by (3) for $b \in B$ that
$$ L^2 \cdot X_b = 0.$$
But for general $b$, the fiber $S_{0,b}$ is an elliptic curve and $X_b$ is a $\bP_1-$bundle over $S_{0,b}$ since $\Delta$ does not meet $S_{0,b}.$
Moreover $L \vert X_b$ is strictly nef, hence ample, contradiction.
\end{proof}

\begin{remark} {\rm  Suppose in (4.2) that $\phi: X \to S$ is a conic bundle, but not
necessarily with $\rho(X/S) = 1.$ Then all arguments still remain valid if $ K_X + t_0 L$ is the $\phi-$pull-back of a $\bQ-$bundle on $S$, for some rational $t_0$. }
\end{remark}

\section{Birational maps}
\setcounter{equation}{0}

In order to prove Conjecture 1.2 in the remaining uniruled cases, it is natural to consider the Mori program. If $X$ admits a contraction
contracting a divisor to a point, the situation is easily understood.

\begin{theorem} Let $X$ be a smooth projective threefold, $L$ strictly nef on $X$.
Suppose that $X$ admits a birational Mori contraction $\phi: X \to Y$ contracting the
exceptional divisor $E$ to a point. Then $K_X  + tL$ is ample for $t > n+1.$
\end{theorem}

\begin{proof} Suppose that $K_X + tL$ is not ample. Write $$ K_X = \phi^*(K_Y) + aE;$$
then $a \in \{2,1,{{1}\over {2}} \}.$ Possibly after replacing $L$ by $2L$ in case $a = {{1}\over {2}},$ we can moreover write
$$ L = \phi^*(L') - bE $$ with a line bundle $L'$
on $Y.$ Notice that $b > 0$ since $L$ is strictly nef. Introduce
$$ D = b K_X + aL; \  D' = bK_Y + aL'.$$
Since $L'$ is again strictly nef, $K_Y + tL'$ is strictly nef for $t >> 0.$ Using (1.6)(1) on $X$ it is a simple matter to verify
$$ (K_Y + tL')^3 > 0$$
for large $t$, so that $K_Y + tL'$ is ample.  Hence we find
positive integers $p,q$ such that $pK_Y + qL'$ is spanned. Choose
$S \in \vert pK_Y  + qL' \vert $ smooth. Now a simple calcluation
shows that
$$ D'^2 \cdot (pK_Y + qL')=  D' \cdot (pK_Y + qL')^2 = 0.$$
Thus $D'_S \cdot (pK_Y + tL')_S = 0.$ Moreover $(D'_S)^2 = 0.$
Hence $D'_S \equiv 0$ by the Hodge index theorem. Thus $D' \equiv
0$ Hence $D \equiv 0$ so that $aL \equiv -bK_X$. Therefore $X$ is
Fano by by Serrano (1.5) and $K_X + tL$ is ample for $t > 4,$
contradiction.
\end{proof}

In case that the contraction $\phi: X \to Y$ contracts a divisor to a curve $C$, the situation is more involved. The reason is that the induced
line bundle $L'$ on $Y$ is not necessarily strictly nef, in fact we can have $L' \cdot C \leq 0.$ We
have already shown that if $X$ admits a Mori fibration or a divisorial
contraction to a point, then the conjecture holds. Since $X$ is
smooth, it remains to consider the case that {\it all}
the extremal rays produce a divisorial contraction to a
nonsingular curve.

\begin{proposition}
Let $X$ be a smooth uniruled threefold, $L$ strictly nef on $X$.
Suppose that all  extremal contractions on $X$ contract a divisor
to a curve. Then $K_X+tL$ is ample for large $t$.
\end{proposition}

\begin{proof} (a)
Let us fix some notations first. Let ${\phi_i},\  i \in I \subset \bN$ be the extremal
contractions on $X,$ with exceptional divisor $E_i$. Let
$C_i:=\phi_i(E_i)$ so that $E_i$ is a
$\mathbb{P}_1$ bundle over $C_i$. Let $[l_i] \in K_X^{<0}$ denotes
the class of the contracted ruling lines in $E_i$.

Let $$\mu:= \min\{ \frac{L \cdot l_i}{-K_X \cdot l_i}\}=\min\{ L
\cdot l_i\} \in \mathbb{N}.$$
Reorder $I$ so that $\phi_1, \ldots, \phi_n$ are exactly those contractions with
$$ L \cdot l_i = \mu.$$

Then the divisor
$$D:= L+ \mu K_X $$
is nef, as a consequence of the cone theorem and the definition of $\mu.$
Moreover, if $D \cdot B = 0$ for some $B \in  \overline{NE}(X),$
then $K_X \cdot B \leq 0.$ In other words,
$$ D^{\perp}
\cap \overline{NE}(X) \subset K_X^{\leq 0}.$$

In particular, if $B$ is an effective curve, then $D \cdot B = 0$ forces
$K_X \cdot B < 0,$ because otherwise $K_X \cdot B = 0,$ hence $L \cdot B = 0,$
contradicting the strict nefness of $L.$
\vskip .2cm \noindent

Our goal is to show that some multiple $mD=mL+m \mu K_X$ is effective, so that we are done by (1.6.3).  \\
Let $\phi = \phi_1 : X \to X_1 = X'$ be the contraction of $E = E_1.$
Let $[l]=[l_1]$ and set $L':= (\phi_*L)^{**}$,
$$D':=L'+\mu K_{X'},   D:=\phi^*(D')=L+\mu K_{X}$$ and let
$C = \phi(E)$.

\vskip .2cm \noindent (b)
We introduce the following numbers
$$\tau := L' \cdot C, \sigma := K_{X'} \cdot C, \gamma = c_1(N^*_{C/X}).$$
Furthermore, let $g$ be the genus of $C$ and $\chi = 2-2g. $

\vskip .2cm \noindent
First we treat the case $L' \cdot C > 0$ so that $L'$ is strictly nef.
Then by induction on $\rho,$ the bundle $K_{X'} + tL'$ is ample, for $t > 4.$
Let $t_0 = {{1} \over {\mu}}.$ Then $K_{X'} + t_0 L'$ is nef, since $D'$ is nef.
Let $\epsilon > 0$ be a small positive number.  Then
$$ K_{X'} + {{t_0} \over {1 - \epsilon}} L'$$
is big (otherwise we would have $(K_{X'} + tL')^3 = 0$ for all $t$ which is absurd).
Now the base point free theorem implies that some multiple $m(K_{X'} + t_0L')$
is spanned, hence $m'D'$ is spanned, and we are done.

\vskip .2cm \noindent

Thus we are reduced to
$$ L' \cdot C \leq 0.$$
Hence $K_{X'} \cdot C \geq 0$, and $C$ is rigid, since $L'.C'>0$ for every irreducible effective curve $C'\neq C$ on $X'$.
We claim that:
$$ D' \cdot C \geq 1. \eqno (*) $$
In fact, we need only to exclude the case: $D' \cdot C = 0.$ Assuming that, we obtain
$$ L' \cdot C + \mu K_{X'} \cdot C = 0$$
and
$$ L_E \equiv - \mu K_X \vert E. $$
Since $L \cdot C_0 > 0,$ we have $K_{X'} \cdot C_0 < 0,$ hence $C_0$ moves. Since $C$ is rigid,
$C_0$ can move only inside $E,$ hence $e \leq 0.$
Write $N^*_E \equiv C_0 + \lambda l.$ Then it is easily checked that
$\lambda = {{1} \over {2}} \gamma + {{1} \over {2}} e $, in the notations of [Ha 77];
so that
$$ N^*_E = C_0 + ({{1} \over {2}}\gamma + {{1} \over {2}} e)l. $$
Since $L_E$ is strictly nef, so is $-K_X \vert E - N^*_E = C_0 + (e+2-2g-\lambda)l,$ so that we conclude:
$$ e + 2-2g - {{1} \over {2}}\gamma - {{1} \over {2}} e \geq {{e} \over 2}, \eqno (**) $$
hence $$ 2-2g \geq {{1} \over {2}} \gamma,$$
with strict inequality for $e = 0,$ since on those ruled surfaces all strictly nef line bundles are ample.
\vskip .2cm \noindent
By the adjunction formula we have
$ \gamma = \sigma + (2-2g),$
hence $\sigma \leq 2-2g.$ Since $\sigma \geq 0,$ we obtain $g \leq 1.$ But a strictly nef divisor on a ruled surface over a rational
 or an elliptic curve is ample, hence the inequality (**) is strict. Thus $g = e = 0$ and $\sigma < 2, \gamma \leq 3.$
So
$$ N^*_C = \sO(k) \oplus \sO(k) $$
with $0 < \gamma = 2k \leq 3,$ hence $k = 1$ and $\sigma = 0.$ So $K_{X'} \cdot C = 0 = L' \cdot C,$ and $L'$ is nef.
If for large $t$, the nef bundle $K_{X'} + tL'$ is big, then we conclude as in the case $L' \cdot C > 0.$ So we
may assume that $K_{X'} + tL')^3 = 0$ for all $t$. Then $K_{X'}^3 = 0.$ However $K_X^3 = 0$ forces $K_{X'}^3 = -2,$
contradiction. Thus we must have
$$ D' \cdot C > 0.$$

\vskip .2cm \noindent  (c) Case: $D'^{\perp} \cap {\overline {NE}}(X) \subset  K_{X'}^{\perp}.$ \\
We are going to rule out this case. Assume there is an irreducible curve $B' \in {\overline {NE}}(X')$ such that $D' \cdot B' = 0.$ Necessarily
$B' \ne C.$ By assumption, $K_{X'} \cdot B' = 0.$
Let $B$ be the strict transform of $B'$ in $X.$ Then $D \cdot B = 0.$ Since $E \cdot B \geq 0,$ we also get
$K_X \cdot B \geq 0.$ Since $L \cdot B > 0$ and $D \cdot B = 0,$ this is impossible.
Hence $D'$ is strictly nef and by induction, $K_{X'} + tL'$ is ample for large $t$. On the other hand, $D'$ is not ample, hence there exists
a nonzero class $B^*\in {\overline {NE}}(X')$ with $D' \cdot B^* = 0,$ hence $K_{X'} \cdot B^* = 0,$ by assumption.This is absurd.

\vskip .2cm \noindent (d) Case:  $D'^{\perp} \cap {\overline {NE}}(X) \not\subset  K_{X'}^{\perp}.$ \\
Then we find $B' \in {\overline {NE}}(X')$ such that $D' \cdot B' = 0 $ and $K_{X'} \cdot B' < 0.$ Since $D'$ is nef, we also find an extremal
curve $l'$ with $D' \cdot l' = 0$. Let $\phi': X' \to X''$ be the associated contraction.
\vskip .2cm \noindent
(d.1) Suppose that $\dim X'' \leq 2.$ Observe that $D' = \phi'^*(D'')$ with a nef bundle $D'$ on $X''.$ So if $\dim X'' \leq 1,$ the bundle $D'$ has
a section and we are done. The same argument works if $\dim X'' = 2$ and $D''^2 \ne 0.$ In the remaining case we need more arguments.
Let $l'$ be a smooth conic and assume that $l'$ meets $C.$ Let $l$ be its strict transform in $X.$ Then $K_X \cdot l \geq -1.$ Since
$D \cdot l = 0$ and $L \cdot l > 0,$ necessarily $K_X \cdot l = -1$ and $E \cdot l = 1.$ Thus $l$ meets $C$ transversely in one point.
The same computations show that $C$ cannot meet a singular conic. Thus $C$ is a section of $X' \to X''$ and $X \to X''$ is still a conic bundle.
Then we conclude by Lemma 4.2 and Remark 4.3.
 \vskip .2cm \noindent
(d.2) Suppose $\phi'$ is birational with exceptional divisor $E'.$ \\
If $C \subset  E',$ then, $C$ being rigid, $E'$ must be ruled and $C$ is the exceptional section in $E'.$ Let $l'$ be a ruling line and $l$ its strict
transform in $X.$ Then $K_X \cdot l = 0.$ Since $D \cdot l = 0,$ we have $L \cdot l = 0$, which is absurd. \\
Things are more complicated when $E' \cap C$ is a finite non-empty set. Suppose first that $E'$ is not $\bP_2$ with normal bundle $\sO(-1).$
In this situation we find a rational curve $l' \subset E'$ meeting $C$ with $K_{X'} \cdot l' = -1.$ Let $\hat l$ be the strict transform in $X.$ Then
$$ \phi^*(l') = \hat l + al $$
with some positive integer $a.$ Since $D' \cdot l' = D \cdot l = 0,$ it follows $D \cdot \hat l = 0.$ Now
$$ K_X \cdot \hat l = -1 + a \geq 0.$$
Hence $D \cdot \hat l = (L + \mu K_X) \cdot l > 0,$ contradiction. \\
It remains to do the case $E' = \bP_2$ with normal bundle $\sO(-1). $ Fixing a line $l' \subset E'$ which meets $C,$ the same computations
as above show that $L \cdot \hat l = 1, \mu = 1, K_X \cdot \hat l = -1$ and $a = 1.$ Notice that $E'$ can meet $C$ only in one point (transversely).
In fact, otherwise we choose two points in $E' \cap C$ and a line $l^*$ through these two points. Then the strict transform $\hat l^*$ satisfies
$K_X \cdot \hat l^* \geq  0$, which is impossible, as already observed. Hence $\hat E'$ is ruled over $\bP_1$ with fibers $\hat l.$
Since $\hat E' \cdot \hat l = -1, $ we can blow down $X$
along the projection $\hat E' \to \bP_1$ to obtain $\psi: X \to Y,$ the blow-up of $Y$ along a smooth curve $C' \simeq \bP_1.$
A priori it is not clear that $Y$ is projective. Let $L_Y = (\psi_*(L))^{**}.$ Then
$$ L = \psi^*(L_Y) - \hat E'.$$ Denoting by $C_0$ the exceptional section of $\hat E'$ and noticing that $N^*_{\hat E'} = C_0 + \hat l,$
we obtain
$$ L \vert \hat E' = C_0 + (L_Y \cdot C')+1) \hat l.$$
Since $L \vert \hat E'$ is ample, it follows that $L_Y \cdot C' > 0$ so that $L_Y$ is strictly nef on the Moishezon manifold $Y$. Then $Y$ has
to be
projective: otherwise by [Pe86] we find an irreducible curve $D$ and a positive closed current $T$ on $Y$ such that $[D+T] = 0.$
But $L_Y \cdot D > 0$ and $L_Y \cdot T \geq 0.$
Now, $Y$ being projective, we conclude by the first part of (b).
\\
If finally $E' \cap  C = \emptyset, $ then the strict transform of $E'$ in $X$ is some $E_j, 2 \leq j \leq n,$ hence defines an extremal contraction on $X$ with
the same properties as $\phi$ and we can continue by induction. Since we assume $X$ uniruled, after finitely many steps we arrive
at $\dim X^{[m]} \leq 2$ and argue as above,
\end{proof}

\section{Higher dimensions}
\setcounter{equation}{0}
 \vskip .2cm \noindent In higher
dimensions it is certainly very difficult to deal with Fano
fibrations; however it is instructive to look at $\bP_k-$bundles
to get an idea on the higher dimensional case. Here we can
calculate explicitly.

\begin{theorem} Let $X$ be a $\bP_k-$bundle over a smooth surface
$S.$ Suppose that $L$ is strictly
nef on $X$. Then $K_X+tL$ is ample for $t > k+3.$
\end{theorem}

\begin{proof} After possibly performing a finite \'etale cover, we may
assume that $X$ is the projectivisation of a rank $(r+1)$-bundle
$\sE$ on $S.$ If we allow $\sE$ to be a $\bQ-$bundle, we may
assume that
$$ L = \sO_{\bP(\sE)}(k)$$
with some positive number $k.$
We also introduce $\zeta = \sO_{\bP(\sE)}(1).$
Notice that $\det \sE$ is strictly nef and suppose that $K_X+tL$ is
not ample. Then
$$ K_X^j \cdot L^{r+2-j} = 0$$
for all $j$ by (1.5).
First recall  the following  $$ \zeta^{r+1} - \pi^*c_1(\sE)
\zeta^{r} +\pi^*c_2(\sE) \zeta^{r-1}=0, $$ and $$ K_X= -(r+1)
\zeta + \pi^*(\det \sE + K_S).$$

The equation $L^{r+2}=0$ immediately leads to $$
\zeta^{r+2}=c_1(\sE)^2-c_2(\sE)=0. \eqno(5)$$ Secondly, combining
with $\zeta^{r+2}=0$, the equation $L^{r+1} \cdot K_X=0$ leads to
$$ \zeta^{r+1} \cdot \pi^*(c_1(\sE) + K_S) = c_1(\sE) \cdot  (c_1(\sE) +
K_S)=0.\eqno(6)$$ Moreover, the equation $L^{r} \cdot K_X^2=0$
leads to
$$ \zeta^{r} \cdot \pi^*(c_1(\sE) + K_S)^2 =  (c_1(\sE) +
K_S)^2=0.\eqno(7)$$

By $(6),(7)$, we have $K_S \cdot  (c_1(\sE) + K_S)=0$ and hence
$K_S^2 =c_1(\sE)^2$. Since $\det \sE$ is strictly nef, equation
$(6)$ yields that $K_S^2=c_1(\sE)^2 \ge 0$ and $c_1(\sE) \cdot K_S
\le 0$.

First suppose that $\kappa (S) \geq 0.$ Then  $K_S \cdot \det \sE
= 0$ and $K_S^2 = 0$ for  $\det \sE$ being strictly nef.  Hence
$K_S \equiv 0.$ Then by (1.5) $\det \sE$ is ample, contradicting
$c_1(\sE)^2 = K_S^2.$
\\ It remains to consider $\kappa (S) = - \infty.$ Since $K_S^2 \geq
0,$ $S$ is either rational or a minimal ruled surface over
an elliptic curve. In the latter case, $K_S^2 = 0,$ hence $c_1(\sE)^2
= 0.$ On the other hand, any strictly nef divisor on a ruled surface
over an elliptic curve is ample (use [Ha77,V.2]), a contradiction.
\\ In case of a rational surface $S$, choose a positive integer $m$
such that $m \det \sE$ is Cartier. Then Riemann-Roch and
$(K_S+\det \sE)^2 = 0$ show that $h^0(m(K_S+\det \sE)) > 0.$ This
contradicts via $(K_S+\det \sE) \cdot \det \sE = 0$ the strict
nefness of $\det \sE.$
\end{proof}

\vspace{1cm}
\small
\begin{tabular}{lcl}
\end{tabular}

\vskip .2cm \vskip .2cm \noindent Fr\'ed\'eric Campana,
D\'epartement de Math\'ematiques, Universit\'e de Nancy,  F-54506
Vandoeuvre-les-Nancy, France,\\ frederic.campana@iecn.u-nancy.fr

\vskip .2cm \noindent
Jungkai Alfred Chen,  Department of Mathematics, National Taiwan University, Taipei 106, Taiwan, \\
jkchen@math.ntu.edu.tw

\vskip .2cm \noindent
Thomas Peternell,  Mathematisches Institut,  Universit\" at Bayreuth, D-95440 Bayreuth, \\
thomas.peternell@uni-bayreuth.de

\end{document}